# A Note on the Zero-Free Regions of the Zeta Function
**N. A. Carella, October, 2012**

***Abstract:*** This note contributes a new zero-free region of the zeta function $\zeta(s)$, $s = \sigma + it$, $s \in \mathbb{C}$. The claimed zero-free region has the form $\{\, s \in \mathbb{C} : \Re e(s) = \sigma > \sigma_0 \,\}$, where $\sigma_0 = 21/40$.

## 1. Introduction

This note considers the theory of the zero-free regions of the zeta function $\zeta(s)$. There is a variety of results on the literature on the zero-free regions of the zeta function $\zeta(s)$, $s = \sigma + it$, $s \in \mathbb{C}$. The best-known result claims that the zeta function is zero-free on the complex half plane $\left\{\, s = \sigma + it : \sigma > 1 - c(\log t)^{-2/3}(\log\log t)^{-1/3} \,\right\} \subset \mathbb{C}$, where $c > 0$ is a constant, and $t > t_0$, see [FD], [IK], [IV], [MV] and related literature. A recent result from the theory of primes in short intervals will be employed to derive a new zero-free region of the zeta function $\zeta(s)$ of the form $\{\, s \in \mathbb{C} : \Re e(s) = \sigma > \sigma_0 \,\}$, where $\sigma_0 = 21/40$.

***Theorem* 1.**   Let $s = \sigma + it \in \mathbb{C}$ be a complex number, and define the complex half plane $\mathcal{H}_1 = \{\, s \in \mathbb{C} : \Re e(s) > 21/40 \,\}$. Then $\zeta(s) \neq 0$ for all $s \in \mathcal{H}_1$. In particular, the nontrivial zeros of the zeta function are confined to the subcritical strip $\mathcal{E}_1 = \{\, s \in \mathbb{C} : 19/40 \leq \Re e(s) \leq 21/40 \,\}$.

The proof of this result, which is simple and straightforward, is unfolded in Section 3. Section 4 offers a second proof at the level of the explicit formula. The second proof is shorter than the first proof, but the first proof provides some inside into the distribution of primes in short intervals $(x, x + y]$, $y < x^{1/2}$. In Section 5, the basic idea of an estimate of the Mertens function $M(x) = \sum_{n \leq x} \mu(n)$, and an estimate of its generalization, the Ramanujan function $R(x) = \sum_{n \leq x} c_n(k)$, is put forward and explicated.

## 2. Foundational Materials
This section serves as a pointer to the literature on a few results used in the proof of the new zero-free region. Most of the lengthy and difficult proofs of these results are given in the cited references.

### 2.1 Basic Results On The Zeros Of The Zeta Function
The set of zeros of the zeta function $\zeta(s)$ is defined by $Z = \{\, s \in \mathbb{C} : \zeta(s) = 0 \,\}$. This subset of complex numbers has a disjoint partition as $Z = Z_T \cup Z_N$, the subset of trivial zeros and the subset of nontrivial zeros.

The subset of trivial zeros of the zeta function is completely determined, it is the subset of negative even integers $Z_T = \{ s = -2n : n \geq 1 \text{ and } \zeta(s) = 0 \}$. The trivial zeros are extracted using the functional equation $\Lambda(s) = \Lambda(1-s)$, where $\Lambda(s) = s(s-1)\pi^{-s/2}\Gamma(s/2)\zeta(s)$ for any complex number $s \in \mathbb{C}$, refer to [TT, p. 13], [KA, p. 55], [IV, p. 8]. More precisely,

$$\pi^{-s/2}\Gamma(s/2)\zeta(s) = \pi^{-(1-s)/2}\Gamma((1-s)/2)\zeta(1-s), \tag{1}$$

and the product form of the gamma function is given by

$$\Gamma(s) = \frac{e^{-\gamma s}}{s}\prod_{n \geq 1}\left(1 + \frac{s}{n}\right)^{-1}e^{s/n}, \tag{2}$$

where $\gamma$ is Euler constant, [EL, p. 151]. The gamma function has simple poles at $s = 0, -1, -2, -3, \ldots$. This is well explicated in [RN, p. 222], and other sources. In contrast, the subset of nontrivial zeros

$$Z_N = \{ s \in \mathbb{C} : 0 < \Re e(s) < 1 \text{ and } \zeta(s) = 0 \}$$

remains unknown. It is believed that every nontrivial zero is of the form $s = 1/2 + it$, $t \in \mathbb{R}$. There are formulas for computing these zeros. As of 2012, about 10 trillions zeros has been computed, see [GN], and the vast literature on this subject.

The counting function $N(T)$ tallies the number of zeros of the zeta $\zeta(s)$ function inside the rectangle

$$R(T) = \{ s \in \mathbb{C} : 0 \leq \Re e(s) < 1 \text{ and } 0 \leq \text{Im}(s) \leq T \}$$

on the critical strip $\{ s \in \mathbb{C} : 0 \leq \Re e(s) < 1 \}$.

**Lemma 2.** (vonMangoldt) Let $N(T) = \#\{ s \in R(T) : \zeta(s) = 0 \}$ be the counting function for the number of complex zeros of height $T \geq 0$. Then

$$N(T) = \frac{T}{2\pi}\log\frac{T}{2\pi} - \frac{T}{2\pi} + O(\log T) + S(T),$$

where $S(t) = \pi^{-1}\arg(\zeta(1/2 + it))$, and $T \geq 1$.

For a detailed proof, see [ES, p. 127], [EL, p. 160], [IV, p. 19], and similar literature.

The nonincreasing function $N(\sigma, T) = \#\{ \beta + it : \beta > \sigma,\ 0 < t \leq T,\ \text{and}\ \zeta(\beta + it) = 0 \}$ of the variable $\sigma \geq 1/2$ counts the number of zeros in the rectangle $R(T)$ of a given range of real part $\beta > \sigma$. It satisfies the relations $N(\sigma, T) \leq N(T)$, and $N(\sigma, T) = 0$ for $\sigma \geq 1$. The estimates $N(\sigma, T) \leq AT^{(1-\sigma)B}\log^C T$, where $A$, $B$, $C > 0$ are constants, are called *Zero Density Theorems*, [HL, p. 121], [KA, p. ], [IV, p. 349]. These estimates are used in the theory of primes in short intervals. Id est, to prove the existence of primes in short intervals. Ergo, it is not surprising if the converse results also hold.

**Lemma 3.** The critical zeros $\rho_n = \sigma_n + i\gamma_n$ of the zeta function $\zeta(s)$ satisfy the following.
(i) The imaginary part satisfies $an/\log n \leq \gamma_n \leq bn/\log n$ for some constants $a, b > 0$.
(ii) The imaginary part $\gamma_n \sim 2\pi n/\log n$ as the integer $n \to \infty$.

The proof is a simple consequence of Lemma 2, more details on this result appear in [EL, p. 160], [IV, p. 20], et cetera.

## 2.2. Theta and Psi Functions

Let $x \in \mathbb{R}$ be a real number. The Chebychev theta and psi functions are defined by $\vartheta(x) = \sum_{p \leq x} \log p$ and $\psi(x) = \sum_{p^n \leq x, n \geq 1} n \log p$ respectively. This subsection provides a short introduction to the results on the theta and psi functions. The first result is an estimate of the explicit formula, which is explained in [MO, p. 13], and [RN, p. 244].

**Theorem 4.** Let $x \geq x_0$ be a large real number, and let $T \leq x$. Then

$$\psi(x) = x - \sum_{|\rho| \leq T} \frac{x^\rho}{\rho} + O(xT^{-1}\log^2 x),$$

where $\rho = \sigma + it \in \mathbb{C}$ runs over the nontrivial zeros of the zeta function $\zeta(s)$.

Proof: Confer [KA, p. 69], and similar literature. ∎

Consider the nontrivial zeros $\sigma + it \in \mathbb{C}$ of the zeta function $\zeta(\sigma + it) = 0$, and define the real number

$$\alpha = \sup \{ \sigma : \zeta(\sigma + it) = 0, \ 0 \leq \sigma \leq 1, \text{ and } t \in \mathbb{R} \}.$$

The supremum $\alpha$ of the real part $\sigma$ of the nontrivial zeros of the zeta function on the complex half plane $\{ \Re e(s) = \sigma \geq 0 \}$ is known to satisfy $0 \leq \alpha \leq 1$. Further, by the symmetry of the nontrivial zeros about the critical line $\{ s \in \mathbb{C} : \Re e(s) = 1/2 \}$, it follows that $1/2 \leq \alpha \leq 1$, see [IG, p. 82].

**Theorem 5.** Let $x \geq x_0$ be a large real number. Then

(i) $\psi(x) = x + O(x^\alpha \log x)$,

(ii) $\vartheta(x) = x + O(x^\alpha \log x)$,

(iii) $\psi(x) = x + \Omega_\pm(x^\alpha)$,

(iv) $\pi(x) = li(x) + O(x^\alpha \log x)$.

Proof: These results follow from Theorem 4, and the number $N(T) \leq c_1 T \log T$, $c_1 > 0$ constant, of nontrivial zeros $\rho = \sigma + it$ such that $0 < \sigma < 1$, and $|t| \leq T$, see Lemmas 2 and 3. For further details, confer [IG, p. 83], [KA] and similar literature. ∎

The third part above is due to Grosswald in [GR], it is also discussed in [IV, p. 347]. This is quite similar to other oscillations results on the theta and psi functions such as the following.

**Theorem 6.** For any small real number $\varepsilon > 0$, and a large number $x \geq x_0$, the followings hold.

(i) $\psi(x) = x + \Omega_\pm(x^{\alpha-\varepsilon})$.                  (ii) $\vartheta(x) = x + \Omega_\pm(x^{\alpha-\varepsilon})$.

(iii) $\pi(x) = li(x) + \Omega_\pm(x^{\alpha-\varepsilon})$.

The proof is given in [MV, p. 464]. A special case geared for the assumed supremum $\alpha = 1/2$ is stated below, it also specifies some information on the rate of change of signs.

**Theorem 7.** [KR] There exist a pair of constants $N_0$ and $T_0$ such that for every function $H(T)$ defined for $T \geq T_0$ and satisfying $N_0 \leq H(T) \leq \log T$, there exist a list of points

$2 < x_0 < x_1 < \cdots < x_V \leq T$ , and $V \geq \log T / H(T)$ for which the error term satisfies

$$\left| \psi(x_i) - x_i \right| \geq C x_i^{1/2} \log\log H(T) ,$$

where $C > 0$ is a constant, and the sign function satisfies

$$\text{sgn}(\psi(x_i) - x_i) = -\text{sgn}(\psi(x_{i+1}) - x_{i+1}),$$

for $i = 0, 1, 2, \ldots, V$.

The best unconditional established estimates of the theta and psi functions have subexponentials error terms. These are derived using a zero-free region of the form $\sigma > 1 - c(\log t)^{-2/3}(\log\log t)^{1/3}$, where $c > 0$ is an absolute constant, $t \geq t_0$ is a real number, and the explicit formula, see [KA, p. 227], and [IV, p. 347].

**Theorem 8.** (Prime Number Theorem) Let $x \geq x_0$ be a large number, and let $c > 0$ be an absolute constant. Then

(i) $\psi(x) = x + O(xe^{-c(\log x)^{3/5}(\log\log x)^{-1/5}})$.      (ii) $\vartheta(x) = x + O(xe^{-c(\log x)^{3/5}(\log\log x)^{-1/5}})$.

(iii) $\pi(x) = li(x) + O(xe^{-c(\log x)^{3/5}(\log\log x)^{-1/5}})$.

## 2.3. Primes in Short Intervals

A pair of results that are utilized in the proof of the main result are stated here. These are claims concerning the densities of prime numbers in short intervals.

**Theorem 9.** ([BA]) For all $x \geq x_0$, the interval $(x - x^{525}, x]$ contains prime numbers.

The proof of this result is based on *Sieve Methods*. The analysis of almost all the previous results for primes in short intervals are based on the *Zero Density Theorems*, [IV, p. 269], the explicit formula and related concepts, see [IK], [IV], [KA], and similar references. A list of the improvements made by various workers in the field is complied in [IV, p. 349], for example, the work in [HL, p. 121] uses $N(\sigma, T) \leq AT^{5(1-\sigma)/12} \log^B T$ , $A$, $B > 0$ constants, to achieve $y = x^{7/12+\varepsilon}$ . The earliest result on the existence of primes in the short intervals $(x, x + y]$ of exponential size $y = x^{1-\varepsilon}$, $1/3300 < \varepsilon < 1$, or better, was achieved by Hoheisel. It superseded the previous best result for the existence of primes in the short intervals $(x, x + y]$ of subexponential sizes $y = xe^{-c(\log x)^{1/2-\varepsilon}}$ with $c > 0$ and $\varepsilon > 0$ constants.

**Theorem 10.** ([RE]) There exists an $x_0$ such that for all $x \geq x_0$ and all $y \geq 1$,

$$\pi(x + y) - \pi(x) \leq \frac{2y}{\log y + 3.53}.$$

Combining these results yields the inequality

$$C_0 \frac{y}{\log x} \leq \pi(x + y) - \pi(x) \leq C_1 \frac{y}{\log x},$$

where $C_0, C_1 > 0$ are constants, for the count of primes in the short intervals $(x, x + y]$, $y \geq x^{21/40}$, refer to [HB] for a refinement.

## 3. The Main Result
A zero-free region of the zeta function is a complex half plane of the form

$$\mathcal{H} = \{ s = \sigma + it \in \mathbb{C} : \zeta(s) \neq 0 \text{ for } \mathfrak{Re}(s) = \sigma > 1 - \sigma_0(t) , \ t \in \mathbb{R} \},$$

where $\sigma_0(t) > 0$ is a real valued function or some other definition. The study of the zero-free regions of the zeta function involves elaborate estimates of $\zeta(s)$, $| \zeta(s) |$, its derivative $\zeta'(s)$, $| \zeta'(s) |$, and the exponential sum $\sum_{x < t < 2x} n^{it}$, see [ES, p.88], [KA], [IV, p. 143], and [FD] etc.

The earliest significant zerofree region

$$\mathcal{H}_2 = \{ s \in \mathbb{C} : \zeta(s) \neq 0 \text{ for } \mathfrak{Re}(s) \geq 1 \}$$

is due to delaVallee Poussin, and Hadamard. The proof is given in [IG, p. 28], [JS, p. 106], and similar sources. The most recent zerofree region

$$\mathcal{H}_3 = \{ s \in \mathbb{C} : \zeta(s) \neq 0 \text{ for } \mathfrak{Re}(s) \geq 1 - c(\log t)^{-2/3}(\log\log t)^{-1/3} \}$$

is the accumulation of many years of arduous and extensive works due to various workers, see [FD], [IV], [MV] [MO], et alii.

A different approach explored in this note revolves around a joint application of theta and psi functions and estimates of prime numbers in short intervals. The proof is also based on the symmetry of the two inequalities:

(I) $\quad 0 < y - c_0 x^\alpha \log^2 x \leq \vartheta(x+y) - \vartheta(x) \leq y + c_0 x^\alpha \log^2 x$, $\qquad\qquad$ for $\quad y \geq x^\alpha$, $\qquad\qquad$ (3)

(II) $\quad y - c_0 x^\alpha \log^2 x < 0 \leq \vartheta(x+y) - \vartheta(x) \leq y + c_0 x^\alpha \log^2 x$, $\qquad\qquad$ for $\quad y < x^\alpha$,

where $c_0 > 0$ is a constant. These inequalities split the study of the existence of primes in short intervals into two cases I and II. The symmetry about the supremum $\alpha$ seems to be governed by the functional equation of the zeta function $\Lambda(s) = \Lambda(1-s)$. A closer inspection of inequality II demonstrates why it is difficult to obtain any result on the existence of primes in very short intervals $(x, x+y]$ with $y \leq x^\beta$, $\beta < 1/2$.

Recall that $\alpha = \sup \{ \sigma : \zeta(\sigma + it) = 0, \ 1/2 \leq \sigma \leq 1, \text{ and } t \in \mathbb{R} \}$ is the supremum of the real parts of the nontrivial zeros of the zeta function.

***Theorem*** 1. Let $s = \sigma + it \in \mathbb{C}$ be a complex number, and define the complex half plane $\mathcal{H}_1 = \{ s \in \mathbb{C} : \Re e(s) > 21/40 \}$. Then $\zeta(s) \neq 0$ for all $s \in \mathcal{H}_1$. In particular, the nontrivial zeros of the zeta function are confined to the subcritical strip $\mathcal{E}_1 = \{ s \in \mathbb{C} : 19/40 \leq \Re e(s) \leq 21/40 \}$.

Proof: For a sufficiently large real number $x \geq x_0$, the difference of the theta over the short interval $(x, x+y]$, $y = x^\beta$ with $\beta > 0$, is given by

$$\vartheta(x+y) - \vartheta(x) = x + y + O((x+y)^\alpha \log^2(x+y)) - \left( x + O(x^\alpha \log^2 x) \right)$$
$$= y + O(x^\alpha \log^2 x), \qquad\qquad (4)$$

where the implied constants are distinct, see Theorem 5. By Theorem 9, the interval $(x, x+x^\beta]$ contains primes whenever $\beta \geq 21/40$. This fact immediately yields

$$0 < \vartheta(x+y) - \vartheta(x) = \sum_{x < p \leq x+y} \log p = y + O(x^\alpha \log^2 x). \qquad\qquad (5)$$

In particular,

$$y - c_0 x^\alpha \log^2 x \leq \vartheta(x+y) - \vartheta(x) \leq y + c_0 x^\alpha \log^2 x, \qquad\qquad (6)$$

where $c_0 > 0$ is a constant. Two cases of this inequality will be examined in details.

**Case I.** Let $\alpha \in (1/2, 1]$, let $x \geq x_0$ be a sufficiently large real number, and let $y = x^\beta \geq x^{21/40}$. Consider the inequality

$$0 < y - c_0 x^\alpha \log^2 x \leq \vartheta(x+y) - \vartheta(x) \leq y + c_0 x^\alpha \log^2 x. \qquad\qquad (7)$$

Dividing by $y = x^\beta \geq x^{21/40}$ across the board returns

$$0 < 1 - c_0 x^{\alpha-\beta} \log^2 x$$
$$< \left( \vartheta(x+y) - \vartheta(x) \right) x^{-\beta} \tag{8}$$
$$\leq 1 + c_0 x^{\alpha-\beta} \log^2 x \ .$$

Since the interval $(x, x+y]$ contains prime numbers for any real number $y = x^\beta \geq x^{21/40}$, the left side of inequality (8) implies that $\beta = 21/40 > \alpha$. Ergo, the supremum satisfies the upper estimate

$$\alpha = \sup \ \{ \ \sigma : \zeta(\sigma+it) = 0, \ \ 1/2 \leq \sigma \leq 1, \ \text{and} \ t \in \mathbb{R} \ \} \leq 21/40. \tag{9}$$

Consequently, there are no zeros in the critical substrip

$$\{ \ s \in \mathbb{C} : 21/40 < \Re e(s) \leq 1 \ \}. \tag{10}$$

Simultaneously, the functional equation $\Lambda(s) = \Lambda(1-s)$, see (1), implies that $\{ \ s \in \mathbb{C} : 0 < \Re e(s) \leq 19/40 \ \}$ is also a zerofree region.

**Case II.** Let $\alpha \in (1/2, \ 1]$, let $x \geq x_0$ be a sufficiently large real number, and let $y = x^\beta \geq x^{21/40}$. Consider the reversed inequality

$$y - c_0 x^\alpha \log^2 x < 0 \leq \vartheta(x+y) - \vartheta(x) \leq y + c_0 x^\alpha \log^2 x \ . \tag{11}$$

Dividing by $y = x^\beta \geq x^{21/40}$ across the board returns

$$1 - c_0 x^{\alpha-\beta} \log^2 x < 0$$
$$< \left( \vartheta(x+y) - \vartheta(x) \right) x^{-\beta} \tag{12}$$
$$\leq 1 + c_0 x^{\alpha-\beta} \log^2 x \ .$$

Since $(x, x+y]$ contains prime numbers for any real number $y = x^\beta \geq x^{21/40}$, and the zeta function is zerofree on the region $\{ \ s \in \mathbb{C} : \Re e(s) \geq 1 \ \}$, the left side of inequality (12) implies that the supremum is

$$\alpha = \sup \ \{ \ \sigma : \zeta(\sigma+it) = 0, \ \ 1/2 \leq \sigma \leq 1, \ \text{and} \ t \in \mathbb{R} \ \} = 1. \tag{13}$$

By Theorem 5-iii, the supremum $\alpha = 1$ implies that

$$\psi(x) = x + \Omega_\pm(x), \tag{14}$$

for sufficiently large numbers $x \geq 1$. But this contradicts the Prime Number Theorem, refer to Theorem 8. Specifically, the equations

$$\psi(x) = x + \Omega_\pm(x) \quad \text{and} \quad \psi(x) = x + O(x e^{-c(\log x)^{3/5}(\log \log x)^{-1/5}}) \tag{15}$$

are inconsistent for sufficiently large numbers $x \geq 1$. Further, it renders the Prime Number Theorem as a trivial result since the error term is larger than the main term:

$$\psi(x) = x + O(x^{\alpha} \log\ x) = x + O(x \log^2 x) \,, \tag{16}$$

see the discussion in [EL, p. 168]. Therefore, the correct inequalities are

(I)  $0 < y - c_0 x^{\alpha} \log^2 x \leq \vartheta(x+y) - \vartheta(x) \leq y + c_0 x^{\alpha} \log^2 x \,,$         for   $y \geq x^{\alpha}$, (17)

(II)  $y - c_0 x^{\alpha} \log^2 x < 0 \leq \vartheta(x+y) - \vartheta(x) \leq y + c_0 x^{\alpha} \log^2 x \,,$         for   $y < x^{\alpha}$,

for some $\alpha \in [1/2, 21/40]$.                                                                    ■

In synopsis, the nontrivial zeros of the analytic continuation

$$\zeta(s) = (1 - 2^{1-s})^{-1} \sum_{n \geq 1} (-1)^{n-1} n^{-s} \,, \tag{18}$$

$\Re e(s) > 0$, of the zeta function $\zeta(s) = \sum_{n \geq 1} n^{-s}$, $\Re e(s) > 1$, are confined to the subcritical strip

$$\mathcal{E}_1 = \{\, s \in \mathbb{C} : 19/40 \leq \Re e(s) \leq 21/40 \,\}. \tag{19}$$

Furthermore, since the zeros of the zeta function are fixed and independent of the value of $x > 1$, a known result on the existence and density of prime numbers in almost all short intervals, see [SL], probably implies that the zeta function has all its nontrivial zeros confined to the subcritical strip

$$\mathcal{E}_s = \{\, s \in \mathbb{C} : 1/2 - \delta \leq \Re e(s) \leq 1/2 + \delta \,\}, \tag{20}$$

$\delta > 0$ is an arbitrarily small number. It is expected that the zeta function has its nontrivial zeros on the critical line $\{\, s \in \mathbb{C} : \Re e(s) = 1/2 \,\}$. The Hadamard product form of the zeta function

$$\zeta(s) = \frac{e^{a+bs}}{s(s-1)\Gamma(s/2)} \prod_{\rho} \left(1 - \frac{s}{\rho}\right)^{-1} e^{s/\rho} \,, \tag{21}$$

where $a = -\log 2\pi$, $b = \log 2\pi - 1 - \gamma/2$, and $\gamma$ is Euler constant, [RN], displays all its zeros in full details.

## 4. Second Proof Of The Main Result

The same result can be extracted at the level of the explicit formula. This approach has certain advantages, and it probably simpler.

**Theorem 1.** Let $s = \sigma + it \in \mathbb{C}$ be a complex number, and define the complex half plane $\mathcal{H}_1 = \{\, s \in \mathbb{C} : \Re e(s) > 21/40 \,\}$. Then $\zeta(s) \neq 0$ for all $s \in \mathcal{H}_1$. In particular, the nontrivial zeros of the zeta function are confined to the subcritical strip $\mathcal{E}_2 = \{\, s \in \mathbb{C} : 19/40 - \varepsilon \leq \Re e(s) \leq 21/40 + \varepsilon \,\}$ for any small number $\varepsilon > 0$.

Proof: For a sufficiently large real number $x \geq x_0$, Theorems 4, and 9 imply that the difference of the psi function over the short interval $(x, x + y]$, $y = x^\beta$ with $\beta \geq 21/40$, satisfies

$$0 \; < \; \frac{\psi(x + y) - \psi(x)}{x^\beta} \; = \; 1 \; - \; \frac{\displaystyle\sum_{\mid \rho \mid \, \leq \, T} \frac{(x + y)^\rho - x^\rho}{\rho}}{x^\beta} \; + \; O(x^{1-\beta} T^{-1} \log^2 x), \tag{22}$$

where $T \geq 1$. Let $\alpha = \sup \{\, \Re e(s) \geq 1/2 : \zeta(s) = 0 \,\} \geq 1/2$ be the supremum. Then, by Theorem 6, the error sum satisfies the omega relation

$$\sum_{\mid \rho \mid \, \leq \, T} \frac{(x + y)^\rho - x^\rho}{\rho} = \Omega_{\pm}(x^{\alpha - \varepsilon}), \tag{23}$$

for any small number $\varepsilon > 0$. But, since the left side of (22) is nonnegative, it forces the inequality

$$\frac{\displaystyle\sum_{\mid \rho \mid \, \leq \, T} \frac{(x + y)^\rho - x^\rho}{\rho}}{x^\beta} = \frac{\Omega_{\pm}(x^{\alpha - \varepsilon})}{x^\beta} < 1, \tag{24}$$

for $T \geq x^{1/2} \geq 1$. This immediately leads to the exponents relation $\alpha - \varepsilon - \beta \leq 0$. Therefore, $\alpha \leq 21/40 + \varepsilon$. Consequently, the complex half plane

$$\mathcal{H}_1 = \{\, s \in \mathbb{C} : \Re e(s) > 21/40 + \varepsilon \,\} \tag{25}$$

is a zerofree region of the zeta function. In addition, the functional equation $\Lambda(s) = \Lambda(1 - s)$ implies that the nontrivial zeros of the zeta function are confined to the critical substrip

$$\mathcal{E}_2 = \{\, s \in \mathbb{C} : 19/40 - \varepsilon \leq \Re e(s) \leq 21/40 + \varepsilon \,\} \tag{26}$$

for any small number $\varepsilon > 0$. ∎

## 5. Applications

By means of the explicit formula, the zero-free region $\{\,\sigma > 1 - c(\log t)^{-\kappa}\,\}$, $\kappa > 0$, is mapped to an approximation of the theta function of the shape $\vartheta(x) = x + O(xe^{-c(\log x)^{1/(1+\kappa)}})$, see [IG, p. 56], [IK, p. 60]. In contrast, the zero-free region $\{\,\sigma > \sigma_0\,\}$, $\sigma_0 \in \mathbb{R}$, has a simpler correspondence, namely, $\vartheta(x) = x + O(x^{\sigma_0})$.

**Corollary 11.** Let $x \geq x_0$ be a large real number. Then

(i) $\psi(x) = x + O(x^{21/40+\varepsilon})$,

(ii) $\vartheta(x) = x + O(x^{21/40+\varepsilon})$,

(iii) $\pi(x) = li(x) + O(x^{21/40+\varepsilon})$,

where $\varepsilon > 0$ is an arbitrarily small real number.

Proof: The first two statements follow from inequality I in (17). The third statement follows from (i) or (ii) and summation by parts. ∎

Both the summatory Mobius function $M(x) = \sum_{n \leq x} \mu(n)$, and the summatory Ramanujan function $R(x) = \sum_{n \leq x} c_n(k)$ are integers valued functions. The later is a generalization of the former. The upper estimates of these finite sums are given below.

**Corollary 12.** Let $k \geq 1$ is a fixed integer, for each $n \in \mathbb{N}$, define the finite sum $c_n(k) = \sum_{\gcd(m,n)=1} e^{i2\pi km/n}$, and let $x \geq x_0$ be a large number. Then

(i) $\displaystyle\sum_{n \leq x} c_n(k) = O(x^{21/40+\varepsilon})$,

(ii) $\displaystyle\sum_{n \leq x} \mu(n) = O(x^{21/40+\varepsilon})$,

where $\varepsilon > 0$ is an arbitrarily small real number.

Proof (i): Let $R(x) = \sum_{n \leq x} c_n(k)$, and let $\sigma_s(n) = \sum_{d \mid n} d^s$, where $k \geq 1$ is a fixed integer. Consider the Ramanujan expansion

$$f(s) = \sum_{n \geq 1} \frac{c_n(k)}{n^s} = \frac{\sigma_{1-s}(k)}{\zeta(s)} \tag{27}$$

of a complex variable $s \in \mathbb{C}$, $\Re e(s) > 0$, [EL, p. 75]. By Theorem 1, the analytic continuation of the zeta function is zerofree in the complex half plane $\{\, s \in \mathbb{C} : \Re e(s) > 21/40\,\}$. Thus, $f(s)$ is an analytic and bounded function of the complex variable $s \in \{\, s \in \mathbb{C} : \Re e(s) > 21/40\,\}$. Moreover, the integral representation

$$\frac{\sigma_{1-s}(k)}{\zeta(s)} = \sum_{n \geq 1} \frac{c_n(k)}{n^s} = \int_1^\infty \frac{dR(t)}{t^s} = \lim_{x \to \infty} \frac{R(x)}{x^s} + s \int_1^\infty \frac{R(t)dt}{t^{s+1}} \qquad (28)$$

is also an analytic and bounded function of the complex variable $s \in \mathbb{C}$, $\Re e(s) > 21/40$. Clearly, this implies that $\sum_{n \leq x} c_n(k) = O(x^{21/40+\varepsilon})$, with $\varepsilon > 0$ an arbitrarily small real number. To settle the second statement (ii), observe that

$$c_n(k) = \sum_{d \,\mid\, \gcd(k,\, n)} \mu(n/d)d \,, \qquad (29)$$

so it is the special case $k = 1$, exempli gratia, $c_n(1) = \mu(n)$ for $n \geq 1$. ∎

## References

[BA] Baker, R. C.; Harman, G.; Pintz, J. The difference between consecutive primes. II. Proc. London Math. Soc. (3) 83 (2001), no. 3, 532-562.

[FD] Ford, Kevin, Vinogradov's integral and bounds for the Riemann zeta function. Proc. London Math. Soc. (3) 85 (2002), no. 3, 565-633.

[ES] Edwards, H. M. Riemann's zeta function. Reprint of the 1974 original, Academic Press, New York; Dover Publications, Inc., Mineola, NY, 2001.

[EL] Ellison, William; Ellison, Fern. Prime numbers. A Wiley-Interscience Publication. John Wiley & Sons, , New York; Hermann, Paris, 1985.

[GN] Xavier Gourdon, The 1013 first zeros of the Riemann Zeta function, and zeros computation at very large height, 2004.

[GR] Grosswald, Emil. Oscillation theorems of arithmetical functions. Trans. Amer. Math. Soc. 126 1967 1–28.

[HB] Heath-Brown, D. R. The number of primes in a short interval. J. Reine Angew. Math. 389 (1988), 22-63.

[HL] Huxley, M. N. On the difference between consecutive primes. Invent. Math. 15 (1972), 164-170.

[HL] Huxley, M. N. The distribution of prime numbers. Large sieves and zero-density theorems. Oxford Mathematical Monographs. Clarendon Press, Oxford, 1972.

[IG] Ingham, A. E. The distribution of prime numbers. Cambridge Tracts in Mathematics and Mathematical Physics, No. 30 Stechert-Hafner, Inc., New York 1990.

[IK] Iwaniec, Henryk; Kowalski, Emmanuel. Analytic number theory. AMS Colloquium Publications, 53. American Mathematical Society, Providence, RI, 2004.

[IV] Ivic, Aleksandar, The Riemann zeta-function. Theory and applications. Wiley, New York; Dover Publications, Inc., Mineola, NY, 2003.

[JS] Jameson, G. J. O. The prime number theorem. London Mathematical Society Student Texts, 53. Cambridge University Press, Cambridge, 2003.

[KA] Karatsuba, Anatolij A. Basic analytic number theory. Translated from the second (1983) edition. Springer-Verlag, Berlin, 1993.

[KR] Kaczorowski, Jerzy; Wiertelak, Kazimierz Oscillations of a given size of some arithmetic error terms. Trans. Amer. Math. Soc. 361 (2009), no. 9, 5023–5039.

[MO] Moreno, Carlos Julio Advanced analytic number theory: L-functions. Mathematical Surveys and Monographs, 115. American Mathematical Society, Providence, RI, 2005.

[MV] Montgomery, Hugh L.; Vaughan, Robert C. Multiplicative number theory. I. Classical theory. Cambridge University Press, Cambridge, 2007.

[RE] Ramaré, Olivier; Schlage-Puchta, Jan-Christoph. Improving on the Brun-Titchmarsh theorem. Acta Arith. 131 (2008), no. 4, 351-366.

[RN] Ribenboim, Paulo, The new book of prime number records, Berlin, New York: Springer-Verlag, 1996.

[SL] Selberg, Atle, On the normal density of primes in small intervals, and the difference between consecutive primes. Arch. Math. Naturvid. 47, (1943). no. 6, 87-105.

[TT] Titchmarsh, E. C. The theory of the Riemann zeta-function. Second edition. Edited and with a preface by D. R. Heath-Brown. The Clarendon Press, Oxford University Press, New York, 1986.